\documentclass{article}%
\usepackage{graphicx}
\usepackage[intlimits]{amsmath}
\usepackage{latexsym}
\usepackage{amsfonts}
\usepackage{amssymb}%
\makeatletter
\newfont{\footsc}{cmcsc10 at 8truept}
\newfont{\footbf}{cmbx10 at 8truept}
\newfont{\footrm}{cmr10 at 10truept}
\newtheorem{theorem}{Theorem}

\newtheorem{conjecture}[theorem]{Conjecture}
\newtheorem{corollary}[theorem]{Corollary}

\newtheorem{lemma}[theorem]{Lemma}

\newtheorem{proposition}[theorem]{Proposition}

\newenvironment{proof}[1][Proof]{\noindent{\textbf {#1}  }}  {\hfill$\Box$\bigskip}

\begin{document}

\title{An asymptotically tight bound on the $Q$-index of graphs with forbidden cycles}
\author{Vladimir Nikiforov\thanks{Department of Mathematical Sciences, University of
Memphis, Memphis TN 38152, USA; \textit{email: vnikifrv@memphis.edu}}}
\maketitle

\begin{abstract}
Let $G$ be a graph of order $n$ and let $q\left(  G\right)  $ be that largest
eigenvalue of the signless Laplacian of $G.$ In this note it is shown that if
$k\geq2$, $n>5k^{2},$ and $q\left(  G\right)  \geq n+2k-2,$ then $G$ contains
a cycle of length $l$ for each $l\in\left\{  3,4,\ldots,2k+2\right\}  .$ This
bound on $q\left(  G\right)  $ is asymptotically tight, as the graph
$K_{k}\vee\overline{K}_{n-k}$ contains no cycles longer than $2k$ and
\[
q\left(  K_{k}\vee\overline{K}_{n-k}\right)  >n+2k-2-\frac{2k\left(
k-1\right)  }{n+2k-3}.
\]
The main result of this note gives an asymptotic solution to a recent
conjecture about the maximum $q\left(  G\right)  $ of a graph $G$ with
forbidden cycles. The proof of the main result and the tools used therein
could serve as a guidance to the proof of the full conjecture.\medskip



\end{abstract}

\section{Introduction}

Given a graph $G,$ the $Q$-index of $G$ is the largest eigenvalue $q\left(
G\right)  $ of its signless Laplacian $Q\left(  G\right)  $. In this note we
give an asymptotically tight upper bound on $q\left(  G\right)  $ of a graph
$G$ of a given order, with no cycle of specified length. Let us start by
recalling a general problem in spectral extremal graph theory:\medskip

\emph{How large can }$q\left(  G\right)  $ \emph{be if }$G$\emph{ is a graph
of order }$n,$\emph{ with no subgraph isomorphic to some forbidden graph
}$F?\medskip$

This problem has been solved for several classes of forbidden subgraphs; in
particular, in \cite{FNP13} it has been solved for forbidden cycles $C_{4}$
and $C_{5}.$ In addition, it seems a folklore result that $q\left(  G\right)
>n$ implies the existence of $C_{3},$ and this bound is exact in view of the
star of order $n$. For longer cycles, a general conjecture has been stated in
\cite{FNP13}, which we reiterate next to clarify the contribution of the
present note.

Let $S_{n,k}$ be the graph obtained by joining each vertex of a complete graph
of order $k$ to each vertex of an independent set of order $n-k;$ in other
words, $S_{n,k}=K_{k}\vee\overline{K}_{n-k}.$ Also, let $S_{n,k}^{+}$ be the
graph obtained by adding an edge to $S_{n,k}.$

\begin{conjecture}
\label{con1} Let $k\geq2$ and let $G$ be a graph of sufficiently large order
$n.$ If $G$ has no $C_{2k+1},$ then $q\left(  G\right)  <q\left(
S_{n,k}\right)  ,$ unless $G=S_{n,k}.$ If $G$ has no $C_{2k+2},$ then
$q\left(  G\right)  <q\left(  S_{n,k}^{+}\right)  ,$ unless $G=S_{n,k}^{+}.$
\end{conjecture}

Conjecture \ref{con1} seems difficult, but not hopeless. It is very likely
that it will be solved completely in the next couple of years. Thus, one of
the goals of this note is to make some suggestions for such a solution and to
emphasize the relevance of some supporting results.

The starting point of our work is the observation that both $q\left(
S_{n,k}\right)  $ and $q\left(  S_{n,k}^{+}\right)  $ are very close to
$n+2k-2$ whenever $n$ is large. In fact, the difference between these values
is $\Omega\left(  1/n\right)  ,$ as can be seen from the following proposition.

\begin{proposition}
\label{pro1}If $k\geq2$ and $n>5k^{2},$ then
\[
n+2k-2-\frac{2\left(  k^{2}-k\right)  }{n+2k+2}>q\left(  S_{n,k}^{+}\right)
>q\left(  S_{n,k}\right)  >n+2k-2-\frac{2\left(  k^{2}-k\right)  }{n+2k-3}.
\]

\end{proposition}

These bounds prompt a weaker, yet asymptotically tight version of Conjecture
\ref{con1}, which we shall prove in this note.

\begin{theorem}
\label{th1}Let $k\geq2,$ $n>6k^{2},$ and let $G$ be a graph of order $n.$ If
$q\left(  G\right)  \geq n+2k-2,$ then $G$ contains cycles of length $2k+1$
and $2k+2.$
\end{theorem}

Before going further, let us note a corollary of Theorem \ref{th1}.

\begin{corollary}
Let $k\geq2,$ $n>6k^{2},$ and let $G$ be a graph of order $n.$ If $q\left(
G\right)  \geq n+2k-2,$ then $G$ contains a cycle of length $l$ for each
$l\in\left\{  3,4,\ldots,2k+2\right\}  .$
\end{corollary}

Indeed, if $l\geq5,$ the conclusion follows immediately from Theorem
\ref{th1}. For $l\in\left\{  3,4\right\}  ,$ recall the bound
\[
q\left(  G\right)  \leq\max\left\{  d_{u}+d_{v}:\left\{  u,v\right\}  \in
E\left(  G\right)  \right\}  .
\]
In view of $q\left(  G\right)  \geq n+2,$ there must be an edge $\left\{
u,v\right\}  $ belonging to two triangles; hence $G$ contains both $C_{3}$ and
$C_{4}.\medskip$

Even though Theorem \ref{th1} is weaker than Conjecture \ref{con1}, our proof
is not too short. To emphasize its structure, we have extracted a few
important points into separate statements, which we give next.

\begin{lemma}
\label{le1}If $G$ is a graph with no $P_{2k+1},$ then for each component $H$
of $G,$ either $v\left(  H\right)  =2k$ or $e\left(  H\right)  \leq\left(
k-1\right)  v\left(  H\right)  .$
\end{lemma}

Write $K_{2k}+v$ for the graph obtained by joining a vertex $v$ to a single
vertex of the complete graph $K_{2k}$.

\begin{lemma}
\label{le2}Let $v$ be a vertex of a graph $G$ of order $n.$ If $G$ contains no
$P_{2k+1}$ with both endvertices different from $v,$ then%
\[
2e\left(  G\right)  -d_{v}\leq\left(  2k-1\right)  \left(  n-1\right)  ,
\]
unless $G$ is a union of several copies of $K_{2k}$ and one $K_{2k}+v.$
\end{lemma}

We also need bounds on $q\left(  G\right)  $ for some special classes of
graphs. Since the known upper bounds did not work in these cases, we came up
with a few technical results giving the required bounds.

\begin{lemma}
\label{le3}Let the integers $k,p,m,$ and $n$ satisfy
\[
k\geq2,\text{ \ }m\geq1,\text{ \ }p\geq0,\text{ \ }n=2kp+m,\text{ \ \ }%
n\geq6k+13.
\]
Let $H$ be a graph of order $m$ and let $F$ be the union of $p$ disjoint
graphs of order $2k,$ which are also disjoint from $H$. Let $G$ be the graph
obtained by taking $F\cup H$ and joining some vertices of $F$ to a single
vertex $w$ of $H$. If%
\begin{equation}
q\left(  H\right)  \leq m+2k-2+\frac{6pk}{n+3}, \label{co1}%
\end{equation}
then $q\left(  G\right)  \leq n+2k-2,$ with equality holding if and only if
equality holds in (\ref{co1}).
\end{lemma}

The reason for Lemma \ref{le3} being so technical is that it must support the
proof of the following two quite different corollaries.

\begin{corollary}
\label{cor1}Let $k,p,$ and $n$ be integers such that $k\geq2$ and $n=2\left(
p+1\right)  k+2.$ Let
\[
G=K_{1}\vee\left(  \left(  pK_{2k}\right)  \cup K_{2k+1}\right)  .
\]
If $n\geq6k+13,$ then $q\left(  G\right)  <n+2k-2.$
\end{corollary}

Given a graph $G$ and $u\in V\left(  G\right)  ,$ write $G-u$ for the graph
obtained by removing the vertex $u.$

\begin{corollary}
\label{cor2}Let $k\geq2,$ $G$ be a graph of order $n,$ and $w\in V\left(
G\right)  .$ Suppose that for each component $C$ of $G-w,$ either $v\left(
C\right)  =2k$ or $e\left(  C\right)  \leq\left(  k-1\right)  v\left(
C\right)  .$ If $n\geq6k+13,$ then $q\left(  G\right)  <n+2k-2.$
\end{corollary}

\medskip

In the next section we outline some notation and results needed in our proofs.
The proofs themselves are given in Section \ref{pf}.

\section{Notation and supporting results}

For graph notation and concepts undefined here, we refer the reader to
\cite{Bol98}. For introductory and reference material on the signless
Laplacian see the survey of Cvetkovi\'{c} \cite{C10} and its references. In
particular, let $G$ be a graph, and $X$ and $Y$ be disjoint sets of vertices
of $G.$ We write:

- $V\left(  G\right)  $ for the set of vertices of $G,$ $E\left(  G\right)  $
for the set of edges of $G,$ and $e\left(  G\right)  $ for $\left\vert
E\left(  G\right)  \right\vert $;

- $G\left[  X\right]  $ for the graph induced by $X,$ and $e\left(  X\right)
$ for $e\left(  G\left[  X\right]  \right)  ;$

- $e\left(  X,Y\right)  $ for the number of edges joining vertices in $X$ to
vertices in $Y;$

- $\Gamma_{u}$ for the set of neighbors of a vertex $u,$ and $d_{u}$ for
$\left\vert \Gamma_{u}\right\vert .\medskip$

We write $P_{k},$ $C_{k},$ and $K_{k}$ for the path, cycle, and complete graph
of order $k.\medskip$

Given a graph $G$ and a vertex $u\in V\left(  G\right)  ,$ note that
\[
\sum_{v\in\Gamma_{u}}d_{v}=2e\left(  \Gamma_{u}\right)  +e\left(  \Gamma
_{u},V\left(  G\right)  \backslash\Gamma_{u}\right)  .
\]
Below we shall use this fact without reference.\bigskip

\subsection{Some useful theorems}

Here we state several known results, all of which are used in the proof of
Theorem \ref{th1}. We start with two classical theorems of Erd\H{o}s and
Gallai \cite{ErGa59}.

\begin{theorem}
\label{EGp}Let $k\geq1.$ If $G$ is a graph of order $n,$ with no $P_{k+2},$
then $e\left(  G\right)  \leq kn/2,$ with equality holding if and only if $G$
is a union of disjoint copies of $K_{k+1}.$
\end{theorem}

\begin{theorem}
\label{EGc}Let $k\geq2.$ If $G$ is a graph of order $n,$ with no $C_{k+1},$
then $e\left(  G\right)  \leq k\left(  n-1\right)  /2,$ with equality holding
if and only if $G$ is a union of copies of $K_{k},$ all sharing a single vertex.
\end{theorem}

For connected graphs Kopylov \cite{Kop77} has enhanced Theorem \ref{EGp} as follows.

\begin{theorem}
\label{Kop} Let $k\geq1,$ and let $G$ be a connected graph of order $n$.

\emph{(i)} If $n\geq2k+2$ and $G$ contains no $P_{2k+2},$ then
\[
e\left(  G\right)  \leq\max\left\{  kn-k\left(  k+1\right)  /2,\binom{2k}%
{2}+\left(  n-2k\right)  \right\}  ;
\]

\emph{(ii)} If $n\geq2k+3$ and $G$ contains no $P_{2k+3},$ then%
\[
e\left(  G\right)  \leq\max\left\{  kn-k\left(  k+1\right)  /2+1,\binom
{2k+1}{2}+\left(  n-2k-1\right)  \right\}  .
\]

\end{theorem}

We refer the reader to the more recent paper \cite{BGLS08}, where the
conditions for equality in Kopylov's bounds are determined as well.

We shall use the following sufficient condition for Hamiltonian cycles, proved
by Ore \cite{Ore60}.

\begin{theorem}
\label{tOB}If $G$ is a graph of order $n\geq3$ and
\[
e\left(  G\right)  >\binom{n-1}{2}+1,
\]
then $G$ has a Hamiltonian cycle.
\end{theorem}

The following structural extension of Theorem \ref{EGp} has been established
in \cite{Nik09}.

\begin{theorem}
\label{Ni} Let $k\geq1$ and let the vertices of a graph $G$ be partitioned
into two sets $A$ and $B$. If
\[
2e\left(  A\right)  +e\left(  A,B\right)  >\left(  2k-1\right)  \left\vert
A\right\vert +k\left\vert B\right\vert ,
\]
then there exists a path of order $2k+1$ with both endvertices in $A.$
\end{theorem}

We finish this subsection with two known upper bounds on $q\left(  G\right)
.$ The proof of Theorem \ref{th1} will be based on a careful analysis of the
following bound on $q\left(  G\right)  ,$ which can be traced back to Merris
\cite{Mer98}. The case of equality was established in \cite{FeYu09}.

\begin{theorem}
\label{tM}For every graph $G,$
\[
q\left(  G\right)  \leq\max\left\{  d_{u}+\frac{1}{d_{u}}\sum_{v\in\Gamma_{u}%
}d_{v}:u\in V\left(  G\right)  \right\}  .
\]
If $G$ is connected, equality holds if and only if $G$ is regular or
semiregular bipartite.
\end{theorem}

Finally, let us mention the following corollary, due to Das \cite{Das04}.

\begin{theorem}
\label{Das}If $G$ is a graph with $n$ vertices and $m$ edges, then
\[
q\left(  G\right)  \leq\frac{2m}{n-1}+n-2,
\]
with equality holding if and only if $G$ is either complete, or is a star, or
is a complete graph with one isolated vertex.
\end{theorem}

\section{\label{pf}Proofs}

In the following proofs there are several instances where the bounds can be
somewhat improved at the price of more involved arguments and calculations.
Such improvements seem not too worthy unless geared towards the complete
solution of Conjecture \ref{con1}. Instead, we tried to keep the exposition
concise, so that the main points are more visible.\medskip

\begin{proof}
[\textbf{Proof of Proposition \ref{pro1}}]It is known that
\[
q\left(  S_{n,k}\right)  =\frac{1}{2}\left(  n+2k-2+\sqrt{\left(
n+2k-2\right)  ^{2}-8\left(  k^{2}-k\right)  }\right)  .
\]
Hence, we see that%
\[
q\left(  S_{n,k}\right)  -\left(  n+2k-2\right)  =-\frac{4\left(
k^{2}-k\right)  }{n+2k-2+\sqrt{\left(  n+2k-2\right)  ^{2}-8\left(
k^{2}-k\right)  }}>-\frac{2\left(  k^{2}-k\right)  }{n+2k-3},
\]
and also
\[
q\left(  S_{n,k}\right)  -\left(  n+2k-2\right)  <-\frac{2\left(
k^{2}-k\right)  }{n+2k-2}.
\]
To bound $q\left(  S_{n,k}^{+}\right)  $ let $\mathbf{x}=\left(  x_{1}%
,\ldots,x_{n}\right)  $ be a unit eigenvector to $q\left(  S_{n,k}^{+}\right)
$ and let $x_{1},\ldots,x_{k\text{ }}$ be the entries corresponding to the
vertices of degree $n-1$ in $S_{n,k}^{+}.$ Let $k+1$ and $k+2$ be the vertices
of the extra edge of $S_{n,k}^{+}.$ By symmetry, $x_{1}=\cdots=x_{k}$ and
$x_{k+1}=x_{k+2}.$ Using the eigenequations for $Q\left(  G\right)  $ and the
fact that
\[
q\left(  S_{n,k}^{+}\right)  >q\left(  S_{n,k}\right)  >n+2k-2-\frac{2\left(
k^{2}-k\right)  }{n+2k-3}>n+k-1,
\]
we see that
\[
x_{k+1}^{2}=\frac{k^{2}x_{1}^{2}}{\left(  q\left(  S_{n,k}^{+}\right)
-k-2\right)  ^{2}}<\frac{k}{\left(  q\left(  S_{n,k}\right)  -k-2\right)
^{2}}<\frac{k}{\left(  n-3\right)  ^{2}}.
\]
On the other hand, comparing the quadratic forms $\left\langle Q\left(
S_{n,k}^{+}\right)  \mathbf{x},\mathbf{x}\right\rangle $ and $\left\langle
Q\left(  S_{n,k}\right)  \mathbf{x},\mathbf{x}\right\rangle $ of the matrices
$Q\left(  S_{n,k}^{+}\right)  $ and $Q\left(  S_{n,k}\right)  $, we see that
\[
q\left(  S_{n,k}^{+}\right)  -\left(  x_{k+1}+x_{k+2}\right)  ^{2}%
=\left\langle Q\left(  S_{n,k}^{+}\right)  \mathbf{x},\mathbf{x}\right\rangle
-\left(  x_{k+1}+x_{k+2}\right)  ^{2}=\left\langle Q\left(  S_{n,k}\right)
\mathbf{x},\mathbf{x}\right\rangle \leq q\left(  S_{n,k}\right)  .
\]
Thus, after some algebra, we get
\begin{align*}
q\left(  S_{n,k}^{+}\right)   &  <q\left(  S_{n,k}\right)  +\frac{4k}{\left(
n-3\right)  ^{2}}<n+2k-2-\frac{2\left(  k^{2}-k\right)  }{n+2k-2}+\frac
{4k}{\left(  n-3\right)  ^{2}}\\
&  <n+2k-2-\frac{2\left(  k^{2}-k\right)  }{n+2k+2},
\end{align*}
completing the proof of Proposition \ref{pro1}.
\end{proof}

\bigskip

\begin{proof}
[\textbf{Proof of Lemma \ref{le1}}]Let $H$ be a component of $G.$ Set
$m=v\left(  H\right)  $ and assume that $m\neq2k.$ We shall show that
$e\left(  H\right)  <\left(  k-1\right)  m.$ If $m\leq2k-1,$ then
\[
e\left(  H\right)  \leq\binom{m}{2}=m\left(  \frac{m-1}{2}\right)  \leq\left(
k-1\right)  m,
\]
as claimed. If $m\geq2k+1$, then clause \emph{(ii) }of\emph{ }Theorem
\ref{Kop} implies that%
\begin{equation}
e\left(  H\right)  \leq\max\left\{  \left(  k-1\right)  m-\left(  \left(
k-1\right)  ^{2}+\left(  k-1\right)  \right)  /2+1,\binom{2k-1}{2}+\left(
m-2k-1\right)  \right\}  . \label{K}%
\end{equation}
This inequality splits into
\[
e\left(  H\right)  \leq\left(  k-1\right)  m-\left(  \left(  k-1\right)
^{2}+\left(  k-1\right)  \right)  /2+1\leq\left(  k-1\right)  m,
\]
and
\[
e\left(  H\right)  \leq\binom{2k-1}{2}+\left(  m-2k+1\right)  =\left(
2k-1\right)  \left(  k-2\right)  +m\leq\left(  k-1\right)  m.
\]
Thus, in all cases we see that $e\left(  H\right)  \leq\left(  k-1\right)  m,$
completing the proof of Lemma \ref{le1}.
\end{proof}

\bigskip

\begin{proof}
[\textbf{Proof of Lemma \ref{le2}}]Assume for a contradiction that
\[
2e\left(  G\right)  -d_{v}\geq\left(  2k-1\right)  \left(  n-1\right)  +1,
\]
and that $G$ has no path of order $2k+1$ with both endvertices different from
$v.$ Write $H$ for the component containing $v$ and let $F$ be the union of
the other components of $G$. Since $P_{2k+1}\nsubseteq F$, Theorem \ref{EGp}
implies that
\begin{equation}
2e\left(  F\right)  \leq\left(  2k-1\right)  v\left(  F\right)  , \label{in1}%
\end{equation}
and so%
\[
2e\left(  H\right)  -d_{v}\geq\left(  2k-1\right)  \left(  v\left(  H\right)
-1\right)  +1.
\]
Noting that
\[
2e\left(  H\right)  -d_{v}=\sum_{u\in V\left(  H\right)  \backslash\left\{
v\right\}  }d_{u}\leq\left(  v\left(  H\right)  -1\right)  ^{2},
\]
we find that $v\left(  H\right)  \geq2k+1.$

Assume that $v\left(  H\right)  \geq2k+2.$ Since
\[
2e\left(  H\right)  \geq\left(  2k-1\right)  \left(  v\left(  H\right)
-1\right)  +1+d_{v}>\left(  2k-1\right)  \left(  v\left(  H\right)  -1\right)
,
\]
Theorem \ref{EGc} implies that $H$ contains a cycle $C$ of order $m\geq2k.$ If
$m\geq2k+1,$ then obviously there is a $P_{2k+1}$ with both endvertices
different from $v,$ so let $m=2k.$ Choose a vertex $w\in V\left(  H\right)  $
such that $w\neq v$ and $w\notin C.$ There exists a shortest path $P$ joining
$w$ to a vertex $u\in C.$ By symmetry, we can index the vertices of $C$ as
$u=u_{1},u_{2},\ldots,u_{2k}.$ Take $u_{0}$ in $P$ at distance $1$ from $C.$
Then the sequences $u_{0},u_{1},u_{2},\ldots,u_{2k}$ and $u_{0},u_{1}%
,u_{2k},\ldots,u_{2}$ induce paths of order $2k+1.$ Since $v$ must be an
endvertex to each of them, we see that $u_{0}=v$. But $w\neq v,$ hence $P$
contains a vertex $u_{-1}$ at distance $2$ from $C.$ Now the sequence
$u_{-1},u_{0},u_{1},u_{2},\ldots,u_{2k-1}$ induces a path of order $2k+1$ with
both endvertices different from $v,$ a contradiction completing the proof
whenever $v\left(  H\right)  \geq2k+2.$

It remains to consider the case $v\left(  H\right)  =2k+1.$ In this case $H$
is not Hamiltonian, as otherwise there is a path of order $2k+1$ with both
endvertices different from $v$; hence, Theorem \ref{tOB} implies that
$e\left(  H\right)  \leq k\left(  2k-1\right)  +1$ and so%
\[
2k\left(  2k-1\right)  +2-d_{v}\geq e\left(  H\right)  -d_{v}\geq\left(
2k-1\right)  2k+1.
\]
This is possible only if $d_{u}=1$ and $e\left(  H\right)  =k\left(
2k-1\right)  +1.$ Since $H-v$ is complete, obviously, $H=K_{2k}+v.$ In
addition, in (\ref{in1}) we have $2e\left(  F\right)  =\left(  2k-1\right)
v\left(  F\right)  ,$ and so the condition for equality in Theorem \ref{EGp}
implies that $G$ is a union of several copies of $K_{2k}$ and one copy of
$K_{2k}+v,$ completing the proof of Lemma \ref{le2}.
\end{proof}

\bigskip

\begin{proof}
[\textbf{Proof of Lemma \ref{le3}}]Let $q:=q\left(  G\right)  ;$ assume for a
contradiction that $q\geq n+2k-2,$ and let $\mathbf{x}=\left(  x_{1}%
,\ldots,x_{n}\right)  $ be a unit eigenvector to $q.$ From the eigenequation
for $Q\left(  G\right)  $ and the vertex $w$ we see that%
\[
\left(  q-n+1\right)  x_{w}\leq\left(  q-d_{w}\right)  x_{w}\leq\sum_{i\in
V\left(  G\right)  \backslash\left\{  w\right\}  }x_{i}\leq\sqrt{\left(
n-1\right)  \left(  1-x_{w}^{2}\right)  },
\]
and in view of $q\geq n+2k-2,$ it follows that%
\begin{equation}
x_{w}^{2}\leq\frac{n-1}{\left(  q-n+1\right)  ^{2}+n-1}<\frac{n-1}{n-1+\left(
2k-1\right)  ^{2}}\leq1-\frac{9}{n+8}. \label{in4}%
\end{equation}

On the other hand, let $u\in V\left(  F\right)  $ be such that $x_{u}%
=\max\left\{  x_{v}:v\in V\left(  F\right)  \right\}  $. Set $x:=x_{u}$ and
note that the eigenequation for $u$ implies that
\[
qx=d_{u}x+\sum_{i\sim u}x_{i}=d_{u}x+x_{w}+\sum_{\left\{  i,u\right\}  \in
E\left(  F\right)  }x_{i}\leq2kx+x_{w}+\left(  2k-1\right)  x=\left(
4k-1\right)  x+x_{w}.
\]
Hence, the inequality $q\geq n+2k-2$ implies that
\[
x\leq\frac{x_{w}}{q-4k-1}\leq\frac{x_{w}}{n-2k-1}.
\]
Next, expanding the quadratic form $\left\langle Q\left(  G\right)
\mathbf{x},\mathbf{x}\right\rangle ,$ we find that
\begin{align*}
q  &  =\sum_{\left\{  i,j\right\}  \in E\left(  G\right)  }\left(  x_{i}%
+x_{j}\right)  ^{2}\leq\sum_{\left\{  i,j\right\}  \in E\left(  G_{0}\right)
}\left(  x_{i}+x_{j}\right)  ^{2}+2kp\left(  x+x_{w}\right)  ^{2}+4p\binom
{2k}{2}x^{2}\\
&  \leq q\left(  G_{0}\right)  +2kp\left(  x+x_{w}\right)  ^{2}+4pk\left(
2k-1\right)  x^{2}\\
&  =q\left(  G_{0}\right)  +2pkx_{w}^{2}+4pkxx_{w}+2pk\left(  4k-1\right)
x^{2}\\
&  \leq q\left(  G_{0}\right)  +2pk\left(  1+\frac{2}{n-2k-1}+\frac
{4k-1}{\left(  n-2k-1\right)  ^{2}}\right)  x_{w}^{2}.
\end{align*}
Now, plugging here the bound (\ref{in4}), we get%
\begin{align}
q  &  \leq q\left(  G_{0}\right)  +2pk\left(  1+\frac{3}{n-2k-1}\right)
\left(  1-\frac{9}{n+8}\right) \nonumber\\
&  \leq q\left(  G_{0}\right)  +2pk+6pk\left(  \frac{1}{n-2k-1}-\frac{3}%
{n+8}\right) \nonumber\\
&  =q\left(  G_{0}\right)  +2pk-6pk\left(  \frac{2n-6k-13}{\left(
n-2k-1\right)  \left(  n+8\right)  }\right)  . \label{in6}%
\end{align}
Note that, in view of $n\geq6k+13$ and $k\geq2,$ we have
\[
\frac{2n-6k-13}{\left(  n-2k-1\right)  \left(  n+8\right)  }\geq\frac
{n}{\left(  n-2k-1\right)  \left(  n+8\right)  }\geq\frac{n}{\left(
n-5\right)  \left(  n+8\right)  }>\frac{1}{n+3}.
\]
Plugging this inequality back in (\ref{in6}) and using (\ref{co1}), we obtain
\begin{align*}
n+2k-2  &  \leq q\leq q\left(  G_{0}\right)  +2pk-\frac{6pk}{n+3}\leq
m+2k-2+\frac{6pk}{n+3}+2pk-\frac{6pk}{n+3}\\
&  =n+2k-2.
\end{align*}
Hence $q\leq n+2k-2,$ with equality holding if and only if equality holds in
(\ref{co1}). The proof of Lemma \ref{le3} is completed.
\end{proof}

\bigskip

\begin{proof}
[\textbf{Proof of Corollary \ref{cor1}}]We shall apply Lemma \ref{le3} with
$H=K_{2k+2}$ and $F=pK_{2k}.$ Clearly $2pk=n-2k-2$ and so%
\begin{align*}
q\left(  H\right)   &  =q\left(  K_{2k+2}\right)  =4k+2<v\left(  H\right)
+2k-2+\frac{3\left(  n-2k-2\right)  }{n+3}\\
&  =v\left(  H\right)  +2k-2+\frac{6kp}{n+3}.
\end{align*}
In the derivation above we use that the inequality $n\geq6k+13$ implies that
$3n-6k-6>2\left(  n+3\right)  .$ The conditions for Lemma \ref{le3} are met
and so $q\left(  G\right)  <n+2k-2,$ completing the proof of Corollary
\ref{cor1}.
\end{proof}

\bigskip

\begin{proof}
[\textbf{Proof of Corollary \ref{cor2}}]Let $F$ be the union of all components
of $G-w$ having order exactly $2k,$ and let $p$ be their number, possibly
zero. Let $H$ be the graph induced by the vertices in $V\left(  G\right)
\backslash V\left(  F\right)  .$\ Note that the hypothesis of Corollary
\ref{cor2} implies that $e\left(  H-w\right)  \leq\left(  k-1\right)  \left(
m-1\right)  $ and so
\[
e\left(  H\right)  \leq e\left(  H-w\right)  +m-1\leq\left(  k-1\right)
\left(  m-1\right)  +m-1=k\left(  m-1\right)  .
\]
Now, from Theorem \ref{Das} we get
\begin{equation}
q\left(  H\right)  \leq\frac{2e\left(  H\right)  }{m-1}+m-2\leq\frac{2k\left(
m-1\right)  }{m-1}+m-2=v\left(  H\right)  +2k-2\leq v\left(  H\right)
+2k-2+\frac{6kp}{n+3}. \label{i}%
\end{equation}
Since $n\geq6k+13,$ we can apply Lemma \ref{le3}, obtaining
\[
q\left(  G\right)  <n+2k-2,
\]
unless equality holds in (\ref{i}). Equality in (\ref{i}) implies that $p=0,$
that is to say $G=H.$ Also, by the condition for equality in Theorem
\ref{Das}, we see that $G$ is either complete, or is a star, or is a complete
graph with one isolated vertex. Since $q\left(  G\right)  =n+2k-2$, $G$ cannot
be a star. If $G$ is complete, then $n+2k-2=2n-2$ and so $n=2k,$ contradicting
that $n\geq6k+13.$ For the same reason $n+2k-2<2n-4$ and so $G$ cannot be a
complete graph with one isolated vertex either. Corollary \ref{cor2} is proved.
\end{proof}

\bigskip

\begin{proof}
[\textbf{Proof of Theorem \ref{th1}}]For short, set $q:=q\left(  G\right)  $
and $V:=V\left(  G\right)  .$ Assume for a contradiction that $G$ is a graph
of order $n>6k^{2},$ with $q\geq n+2k-2,$ and suppose that $C_{2k+1}\nsubseteq
G$ or $C_{2k+2}\nsubseteq G.$ We may and shall suppose that $G$ is edge
maximal, because edge addition does not decrease the $Q$-index. In particular,
this assumption implies that $G$ is connected.

Let $w$ be a vertex for which the expression%
\[
d_{w}+\frac{1}{d_{w}}\sum_{i\sim w}d_{i}%
\]
is maximal. We shall show that
\begin{equation}
d_{w}+\frac{1}{d_{w}}\sum_{i\sim w}d_{i}\leq n+2k-2. \label{in3}%
\end{equation}
This is enough to prove Theorem \ref{th1}, unless \
\[
q=d_{w}+\frac{1}{d_{w}}\sum_{i\sim w}d_{i}.
\]
However, $G$ is connected, so if equality holds in (\ref{in3}) Theorem
\ref{tM} implies that $G$ is regular or bipartite semiregular; it is not hard
to see that neither of these conditions can hold. Indeed, if $G$ is bipartite,
then $q\leq n.$ If $G$ is regular, then $q=2\delta\leq n,$ as otherwise,
Bondy's theorem \cite{Bol98} implies that $G$ is pancyclic. So to the end of
the proof we shall focus on the proof of (\ref{in3}).

For short, set $A=\Gamma_{w},$ $B=V\left(  G\right)  \backslash\left(
\Gamma_{w}\cup\left\{  w\right\}  \right)  ,$ and $G_{w}=G\left[
V\backslash\left\{  w\right\}  \right]  .$ Obviously, $\left\vert A\right\vert
=d_{w}$ and $\left\vert A\right\vert +\left\vert B\right\vert =n-1.$

First we shall prove that $C_{2k+1}\subset G.$ Assume thus that $C_{2k+1}%
\nsubseteq G;$ clearly $P_{2k}\nsubseteq G\left[  A\right]  ,$ and so Theorem
\ref{EGp} implies that $e\left(  A\right)  \leq\left(  k-1\right)  \left\vert
A\right\vert $. Now
\begin{align*}
d_{w}+\frac{1}{d_{w}}\sum_{i\sim w}d_{i}  &  =\left\vert A\right\vert
+1+\frac{2e\left(  A\right)  +e\left(  A,B\right)  }{\left\vert A\right\vert
}\leq\left\vert A\right\vert +1+\frac{2\left(  k-1\right)  \left\vert
A\right\vert +\left\vert A\right\vert \left\vert B\right\vert }{\left\vert
A\right\vert }\\
&  \leq\left\vert A\right\vert +1+2k-2+\left\vert B\right\vert =n+2k-2.
\end{align*}
This completes the proof that $C_{2k+1}\subset G.$

The proof that $C_{2k+2}\subset G$ is somewhat longer. Assume that
$C_{2k+2}\nsubseteq G$ and note that if $d_{w}\leq2k-1$, then
\[
d_{w}+\frac{1}{d_{w}}\sum_{i\sim w}d_{i}=d_{w}+\Delta\leq2k-1+n-1=n+2k-2,
\]
so (\ref{in3}) holds. Thus, hereafter we shall assume that $d_{w}\geq2k.$

Further, note that the graph $G_{w}$ contains no path with both endvertices in
$A,$ as otherwise $C_{2k+2}\subset G$. Hence, Theorem \ref{Ni} implies that%
\[
2e\left(  A\right)  +e\left(  A,B\right)  \leq\left(  2k-1\right)  \left\vert
A\right\vert +k\left\vert B\right\vert =\left(  k-1\right)  d_{w}+k\left(
n-1\right)  ,
\]
and therefore
\begin{align*}
d_{w}+\frac{1}{d_{w}}\sum_{i\sim w}d_{i}  &  =d_{w}+1+\frac{2e\left(
A\right)  +e\left(  A,B\right)  }{d_{w}}\leq d_{w}+1+\frac{\left(  k-1\right)
d_{w}+k\left(  n-1\right)  }{d_{w}}\\
&  =d_{w}+k+\frac{k\left(  n-1\right)  }{d_{w}}.
\end{align*}
The function $x+k\left(  n-1\right)  /x$ is convex for $x>0;$ hence, the
maximum of the expression
\[
d_{w}+\frac{k\left(  n-1\right)  }{d_{w}}%
\]
is attained for the minimum and maximum admissible values for $d_{w}.$ Since
$d_{w}\geq2k,$ in either case we find that
\[
d_{w}+\frac{1}{d_{w}}\sum_{i\sim w}d_{i}<n+2k-2,
\]
unless $d_{w}\geq n-2.$ Therefore, to complete the proof we only need to
consider the cases $d_{w}=n-2$ and $d_{w}=n-1.$

First, suppose that $d_{w}=n-2$ and let $v$ be the vertex of $G$ such that
$v\neq w$ and $v\notin\Gamma_{w}.$ Note that $G_{w}$ contains no path of order
$2k+1$ with both endvertices different from $v,$ as such a path would make a
$C_{2k+2}$ with $w.$ Therefore, the hypothesis of Lemma \ref{le2} is
satisfied, and so either
\begin{equation}
2e\left(  A\right)  +e\left(  A,B\right)  =2e\left(  G_{w}\right)  -d_{v}%
\leq\left(  2k-1\right)  \left(  n-2\right)  \label{in2}%
\end{equation}
or $G_{w}$ is a union of several copies of $K_{2k}$ and one $K_{2k}+v.$ If
(\ref{in2}) holds, we see that
\begin{align*}
d_{w}+\frac{1}{d_{w}}\sum_{i\sim w}d_{i}  &  \leq n-2+1+\frac{2e\left(
A\right)  +e\left(  A,B\right)  }{n-2}\leq n-1+\frac{\left(  2k-1\right)
\left(  n-2\right)  }{\left(  n-2\right)  }\\
&  =n+2k-2,
\end{align*}
completing the proof of (\ref{in3}). On the other hand, if $G_{w}$ is a union
of several copies of $K_{2k}$ and one $K_{2k}+v,$ then $G$ is a spanning
subgraph of the graph $G^{\prime}=K_{1}\vee\left(  \left(  pK_{2k}\right)
\cup K_{2k+1}\right)  ,$ with $p$ chosen so that $n=2\left(  p+1\right)  k+2.$
Since $n\geq6k^{2}+1\geq6k+13,$ we can apply Corollary \ref{cor1} obtaining
that
\[
q\left(  G\right)  <q\left(  G^{\prime}\right)  <n+2k-2,
\]
which contradicts the assumption and completes the proof of Theorem \ref{th1}
if $d_{w}=n-2$.

Finally, let $d_{w}=n-1.$ Since $G_{w}$ contains no $P_{2k+1},$ Lemma
\ref{le1} implies that for each component $C$ of $G_{w},$ either $v\left(
C\right)  =2k$ or $e\left(  C\right)  \leq\left(  k-1\right)  v\left(
C\right)  .$ Since $n\geq6k^{2}+1\geq6k+13,$ the graph $G$ satisfies the
hypothesis of Corollary \ref{cor2}, and so%
\[
q\left(  G\right)  <n+2k-2,
\]
completing the proof of Theorem \ref{th1}.
\end{proof}

\bigskip

\textbf{Acknowledgement}

Thanks are due to the referee for helpful suggestions.

\end{document}